\documentclass{amsart}
\usepackage{geometry}
\usepackage{epstopdf}
\usepackage[bookmarksnumbered,colorlinks,bookmarks,citecolor=blue,linkcolor=blue,urlcolor=blue,breaklinks,linktocpage]{hyperref}
\usepackage[all]{hypcap}
\usepackage{subfigure}
\usepackage{tikz-cd}
\newcommand{\Z}{\mathbb{Z}}

\theoremstyle{plain}
\newtheorem{theorem}{Theorem}[section]

\theoremstyle{definition}
\newtheorem{definition}[theorem]{Definition}

\theoremstyle{remark}

\newcommand{\m}[1]{\text{\small #1}}
\begin{document}
\title{A group-theoretical classification of three-tone and four-tone harmonic chords}
\author{Jason K.C. Polak}
\date{\today}
\maketitle
\begin{abstract}
We classify three-tone and four-tone chords based on subgroups of the symmetric group acting on chords contained within a twelve-tone scale. The actions are inversion, major-minor duality, and augmented-diminished duality. These actions correspond to elements of symmetric groups, and also correspond directly to intuitive concepts in the harmony theory of music. We produce a graph of how these actions relate different seventh chords that suggests a concept of distance in the theory of harmony.
\end{abstract}
\tableofcontents

\section{Introduction}

Early on in music theory we learn of the harmonic triads: major, minor, augmented, and diminished. Later on we find out about four-note chords such as seventh chords. We wish to describe a classification of these types of chords using the action of the finite symmetric groups.

We represent notes by a number in the set $\Z/12 = \{0, 1,2,\dots,10,11\}$. Under this scheme, for example, $0$ represents $C$, $1$ represents $C\sharp$, $2$ represents $D$, and so on. We consider only pitch classes modulo the octave.

We describe the sounding of simultaneous notes by an ordered increasing list of integers in $\Z/12$ surrounded by parentheses. For example, a major second interval $M2$ would be represented by $(0, 2)$, and a major chord would be represented by $(0,4, 7)$. We denote the set of all simultaneous $k$-notes by $N_k$. So $(0,4,7)\in N_3$.

Each ordered list describing a set of tones gives rise to a partition of $12$. Consider for example the major chord $(0,4,7)$. The differences between these notes is $4, 3, 5$, the last distance of $5$ being the distance from $7$ from $12$. In general, if $(0,a_1,\dots,a_{k-1})$ is a $k$-tone chord then it is associated to a partition of $12$ via the map
\begin{align*}
    (0,a_1,a_2,a_3,\dots,a_{k-1})\mapsto [a_1,a_2-a_1,a_3-a_2,\dots,12-a_{k-1}].
\end{align*}

We denote the set of all such partitions of $12$ by $P_{12}$, and we consider partitions unordered. To write a partition, we use square brackets, so that the partition associated to $(0,4,7)$ is the partition $[3,4,5]$. Since we consider all partitions unordered, we usually write the numbers of a partition in increasing order. Under this scheme, different intervals or chords may correspond to the same partition. For example, the partition $[3,4,5]$ is also associated to $(0, 0+3, 0 +3  +5) = (0,3,8)$.

To summarize, we can say that to each simultaneous playing of notes, we can assign a partition of twelve. Each partition may give rise to more than one way of playing a set of simultaneous notes, but since each way of doing so corresponds to a different way of adding up the partition of twelve, we can apply elements of the appropriate symmetric group to obtain every possible set of simultaneous notes. This is the basis of the classification scheme we introduce. We also use this classification scheme to create a graph showing how each chord is related with regard to each of the introduced group actions.

The method of studying group actions on chords goes back to \cite{riemann1919handbuch}. Our method of relating chords is similar to the one used in \cite{crans2009musical} and, but the authors in that paper are not concerned with pitch classes or classifying from partitions. They also consider triads, whereas the main focus of this work is four-tone harmonic chords or \emph{seventh} chords. The work of \cite{cannas2017group} is focused on seventh chords, but again is more concerned with absolute pitches. Also, unlike in these previous interesting studies, we are more concerned with deriving the basic harmonic three-tone and four-tone chords axiomatically and with representing this relationship graphically.

\section*{Acknowledgements} The author would like to thank Emily Polak for a conversation about naming the augmented-diminished duality operator, and for inspiring this paper with her musical talent.

\section{Three-tone harmonic chords}

As a warmup for three-tone chords, we consider the following definition, which gives us the usual three-tone harmonic chords.
\begin{definition}
    A harmonic three-tone chord is any chord whose partition does not contain any $x \leq 2$.
\end{definition}
Therefore, the major chord $(0,4,7)$, whose partition is $[3,4,5]$, is an example of a harmonic chord since it does not contain $1$ or $2$. We will see now that this definition encompasses the usual harmonic chords that we learn in beginning music theory.

In order to classify harmonic three-tone chords, we first list all partitions of twelve that do not have any $x\leq 2$ in them. There are three of these: $[3, 3,6], [3,4,5],$ and $[4,4,4]$. Each may correspond to many chords; however to list all such chords we just need to find one chord and apply element of the symmetric group $S_3$ to it. We start with $[3,4,5]$. One such corresponding chord is $(0,4,7)$, the major chord. Here are the possibilities we are starting with:
\begin{enumerate}
    \item $(0,4,7) \mapsto [3,4,5]$, the major chord.
    \item $(0,3,6) \mapsto [3,3,6]$, the diminished chord.
    \item $(0,4,8) \mapsto [4,4,4]$, the augmented chord.
\end{enumerate}

We now consider two operations on chords. The first is inversion, named after the corresponding inversion in music theory. It is the function
\begin{align*}
    i_k: N_k&\longrightarrow N_k\\
    (0,a_1,\dots,a_{k-1})&\longmapsto (0, a_2-a_1,\dots, a_{k-1}-a_1, 12-a_1).
\end{align*}
When no confusion is possible, we denote $i_k$ by $i$. Despite its name, inversion does not have order two except when $k=2$, but we keep the terminology to be consistent with music theory. In fact, it is easily seen that $i_k$ has order $k$; that is, $i_k^k = {\rm id}$, where ${\rm id}$ is the identity operator. There is a second operator, which we call major-minor duality, denoted by $d_k$:
\begin{align*}
    d_k: N_k&\longrightarrow N_k\\
    (0,a_1,\dots,a_{k-1})&\longmapsto (0,12-a_{k-1},12-a_{k-2},\dots,12-a_1).
\end{align*}
Again, when $k$ is clear, we denote $d_k$ by $d$. The duality operator satisfies $d^2 = {\rm id}$. In this way, we get an action of the dihedral group $D_k$ on $N_k$. Recall the dihedral group on $k$ vertices is a group of order $2k$ and is given by generators and relations in the presentation
\begin{align*}
    D_k = \langle~r,s ~|~ r^k, s^2, r^ks = sr^{-k}~\rangle.
\end{align*}
The action of $D_k$ on $N_k$ is that of $r$ acting as $i_k$ and $s$ acting as $d_k$.

Let us now return to $k = 3$, the situation of three-tone chords. In this case, $D_3 = S_3$, where in general $S_k$ is the full symmetric group on $k$ symbols. We now recall the chords we listed corresponding to all partitions satisfying the definition of a harmonic three-tone chord. By substituting these chords into the functions of inversion and duality, we will obtain all harmonic three-tone chords by definition. We have for inversion:
\begin{enumerate}
    \item The major chord: under inversion, we obtain the chords
        \begin{align*}
            &\{ (0,4,7), i(0,4,7) = (0,3,8), i(0,3,8) = (0,5,9)\}\\
            =&\{ (0,4,7), (0,3,8), (0,5,9)\}.
        \end{align*}
        by successive application of $i = i_3$. In music theory, these correspond to the root, first inversion, and second inversion respectively of the major chord. Applying the duality operator to the above chords gives:
        \begin{align*}
            &\{ d(0,4,7) = (0,5,8), d(0,3,8) = (0,4,9), d(0,5,9) = (0,3,7)\}\\
            =&\{ (0,5,8), (0,4,9), (0,3,7)\}.
        \end{align*}
        We recognize this as the familiar minor chords and its inversions. In fact, we see from this that through duality, the root, first, and second inversions of the major chord correspond to the second, first, and root positions of the minor chord respectively. This follows from the relation property $di^k = i^{-k}d$ satisfied by the dihedral group. Because the major and minor chords are dual, the $d$ operator could be thought of as major-minor duality.

    \item The diminished chord: under inversion, we obtain the chords:
        \begin{align*}
            \{ (0,3,6), i(0,3,6) = (0,3,9), i(0,3,9) = (0,6,9)\}
        \end{align*}
        by successive applications of $i$. These are the root, first inversion, and second inversion of the diminished chord. Applying the duality operator gives $d(0,3,6) = (0,6,9)$. Therefore, we conclude that the duality operator transposes the root and second inversion of the diminished chord and fixes the first inversion.

    \item The augmented chord: under inversion, we find that
        \begin{align*}
            i(0,4,8) = (0,4,8).
        \end{align*}
        So, under inversion, the augmented chord is stable. In some sense it is this fixed property of the augmented chord that actually demands resolution. Not only that, but $d(0,4,8) = (0,4,8)$. Therefore, the augmented chord is both inversion stable and major-minor dual.
\end{enumerate}

\section{Four-tone harmonic chords}

Four-tone harmonic chords are more complex than three-tone harmonic chords. Recall that we have defined a three-tone harmonic chord to be one whose partition does not contain $1$ or $2$. On the other hand, it makes less sense to define a four-tone harmonic chord in exactly the same way, since that would exclude seventh chords. Instead, we use the following.
\begin{definition}\label{defn:fourtoneharmonic}
    A four-tone harmonic chord is a four-tone chord whose partition contains at most one $x\leq 2$.
\end{definition}

We have already described two operators $i_k$ and $d_k$, called inversion and duality, that will operate on the set $N_k$, the set of $k$-tone chords. However, in the case of $k=4$, the symmetric group $S_4$ has $24$ elements, which is larger than the dihedral group $D_4$ of eight elements. Can we obtain any additional actions on $N_4$, and if so, do those actions contain any new operators of harmonic significance? We will answer this question now.

To begin, we again list all partition of $12$ of length $4$ that satisfy Definition~\ref{defn:fourtoneharmonic} of a four-tone harmonic chord along with one example of each type. To keep the verbiage simple for seventh chords, we describe seventh chords first by their three-tone harmonic chord and the type of seventh added to them as follows:

\begin{enumerate}
 \item $(0,4,7,11)\mapsto [4, 4, 3, 1]$, the major-major seventh (also known as the major seventh chord)
 \item $(0,4,7,10) \mapsto [4, 3, 3, 2]$, the major-minor seventh chord (also known as the dominant seventh chord)
 \item $(0,3,6,9) \mapsto [3, 3, 3, 3]$, the diminished-diminished seventh chord (also known as the fully diminished seventh chord)
\end{enumerate}
We can calculate how many chords we will find before we actually find them. For the permutations $[4,4,3,1]$ and $[4,3,3,2]$ there are $4!/2! = 12$ possible chords, and there is only one possible chord for $[3,3,3,3]$. Therefore, there are $25$ four-tone harmonic chords. As before, we write down all the possible chords we can obtain from these using the inversion and duality operators.

Let us first discuss the major-major seventh. Its inversions are the orbit of the major-major seventh element $(0,4,7,11)$ under the inversion operator. Explicitly:
\begin{align*}
    \{ (0,4,7,11), (0,3,7,8), (0,4,5,9), (0,1,5,8)\}.
\end{align*}
Something else which is quite intriguing is that $d(0,4,7,11) = (0,1,5,8)$. That is $i^3 = d$ when restricted to the orbit of the major-major seventh under $i$. In other words, adding the major seventh interval to the major chord has stabilized the duality between major and minor, so that the dual of the major-major seventh orbit is again the major-major seventh orbit. The harmonic significance is that the inversions of the major-major seventh pair up as follows under duality:
\begin{align*}
    (0,4,7,11)&\leftrightarrow (0,1,5,8)\\
    (0,3,7,8)&\leftrightarrow (0,4,5,9).
\end{align*}
This pairing of inversions means that in harmony, inversions of the major-major seventh according to this pairing also function as major-minor tension as well. This relationship was not present in three-tone chord harmony.

The more prevalent major-minor seventh chord, often functioning as a dominant seventh chord, behaves differently. We still have its inversions:
\begin{align*}
    \{ (0,4,7,10), (0,3,6,8), (0,3,5,9), (0,2,6,9) \}.
\end{align*}
Applying the duality operator to this set gives the following chords and their inversions:
\begin{align*}
    \{ (0,2,5,8), (0,3,6,10), (0,3,7,9), (0,4,6,9)\}.
\end{align*}
The root position of this set (the one which contains either a major or minor seventh interval from the lowest note) is $(0,3,6,10)$, which is the diminished-minor seventh chord. We have thus far only accounted for 9 of the 25 types of four-tone harmonic chords. We recall the reason for this is that the dihedral group $D_4$ only has $8$ elements whereas $S_4$ has $24$ elements. Therefore, we need to write down a new type of operator on chords to obtain the remaining types of four-tone harmonic chords.

We consider a new operator on four-tone chords called augmented-diminished duality:
\begin{align*}
    a:N_4&\longrightarrow N_4\\
    (0,a_1,a_2,a_3)&\longmapsto (0,a_1,a_1 + a_3 - a_2,a_3).
\end{align*}
The reader will verify that $a^2 = {\rm id}$, similar to duality. We note that we can also specify chords of the form $(0,a_1,a_2,a_3)$ by an \emph{ordered partition}. That is, the chord $(0,a_1,a_2,a_3)$ corresponds to the \emph{ordered} partition $[a_1,a_2-a_1,a_3-a_2,12-a_3]$. The augmented-diminished operator $a$ then corresponds to the map
\begin{align*}
    [a,b,c,d]\mapsto [a,c,b,d].
\end{align*}
The augmented-diminished map is so-named because it transforms the major-major seventh inversions as follows:
\begin{align*}
    \{ (0,4,7,11), (0,3,7,8), (0,4,5,9), (0,1,5,8)\}\\
    \downarrow\\
    \{ (0,4,8,11), (0,3,4,8), (0,4,8,9), (0,1,4,8)\}.   
\end{align*}
That is, the augmented-diminished operator $a$ satisfies $a(0,4,7,11) = (0,4,8,11)$. It transforms the major-major seventh into a augmented-major seventh. We note that the second set is not a set of inversions, although it contains some inversion pairs. Namely, the second inversion of $(0,4,8,11)$ is $(0,3,4,8)$. What about the chord $(0,4,8,9)$? The inversions of this chord are:
\begin{align*}
    \{ (0,3,7,11), (0,4,8,9), (0,4,5,8), (0,1,4,8)\}.
\end{align*}
We see two chords that were in the augmentation of the inversions of the major-major seventh. The chord $(0,3,7,11)$ is a new seventh chord in our list, the minor-major seventh chord. On the other hand, the set of inversions of the augmented-major seventh chord are
\begin{align*}
    \{ (0,4,8,11), (0, 4,7,8), (0,3,4,8), (0,1,5,9)\}.
\end{align*}
Finally, if we take the diminished-minor chord $(0,3,6,10)$ and apply the augmented-diminished operator we get $(0,3,7,10)$, a minor-minor chord. Its inversions are
\begin{align*}
    \{ (0,3,7,10), (0,4,7,9), (0,3, 5, 8), (0,2,5,9)\}.
\end{align*}
This finishes our analysis of all four-tone harmonic chords. We have summarized these classified into inversion orbits in Table~\ref{tab:fourtoneInversionOrbit}.

\begin{table}[h]
    \centering
    \begin{tabular}{|l | l|}
        \hline
        Chord type & Inversions \\
        \hline\hline
        Major-Major & $\{ (0,4,7,11), (0,3,7,8), (0,4,5,9), (0,1,5,8)\}$ \\
        \hline
        Minor-Major & $\{ (0,3,7,11), (0,4,8,9), (0,4,5,8), (0,1,4,8)\}$\\
        Augmented-Major & $\{ (0,4,8,11), (0, 4,7,8), (0,3,4,8), (0,1,5,9)\}$\\
        \hline\hline
        Major-Minor & $\{ (0,4,7,10), (0,3,6,8), (0,3,5,9), (0,2,6,9) \}$\\
        Diminished-Minor & $\{ (0,3,6,10), (0,3,7,9), (0,4,6,9), (0,2,5,8)\}$\\
        \hline
        Minor-Minor & $\{ (0,3,7,10), (0,4,7,9), (0,3, 5, 8), (0,2,5,9)\}$\\
        \hline\hline
        Diminished-Diminished & $\{ (0,3,6,9) \}$\\
        \hline
    \end{tabular}
    \caption{A table of all four-tone harmonic chords, where each row corresponds to an orbit under the inversion map. The single horizontal lines in the table group together dual sets of inversions; that is inversion that are taken to each other under major-minor duality. The double line separates the three orbits of the augmented-diminished operator, which also correspond to the three permutations that we originally used to derive all these chords.}
    \label{tab:fourtoneInversionOrbit}
\end{table}

\section{The chord graph}

Recall that we have used the augmented-diminished operator
\begin{align*}
    a:N_4&\longrightarrow N_4\\
    (0,a_1,a_2,a_3)&\longmapsto (0,a_1,a_1 + a_3 - a_2,a_3).
\end{align*}
along with inversion and duality to derive all four-tone harmonic chords from the three permutations $[4,4,3,1], [4,3,3,2],$ and $[3,3,3,3]$. In this section we talk a little about the harmonic significance of the augmented-diminished operator. In order to visualize more readily how the augmented-diminished operator acts, we examine the following directed graph showing the relations between all four-tone harmonic chords, where the relations are defined by the operators inversion, major-minor duality, and augmented-diminished duality:\\
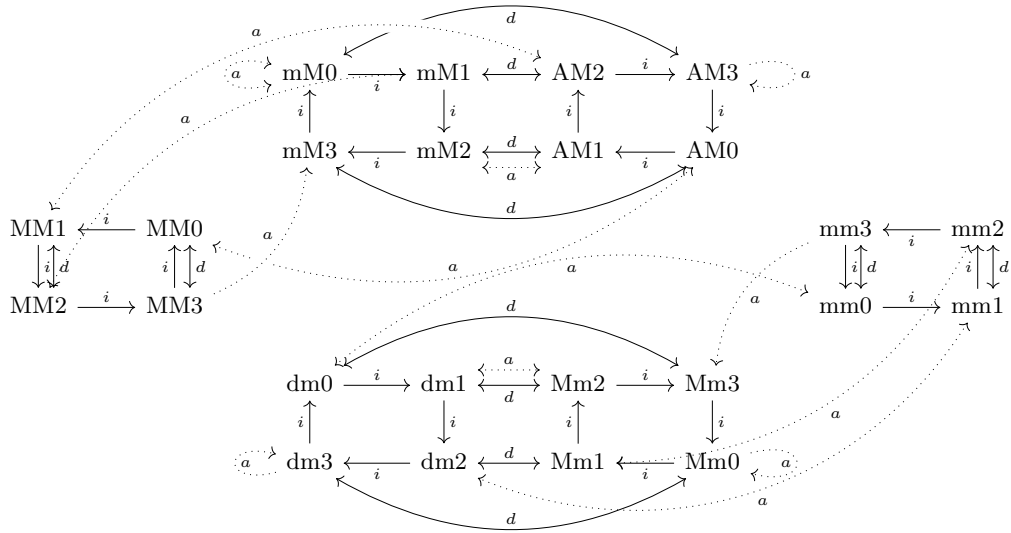
\begin{center}
{\small
\begin{tikzcd}
    ~ & ~ & \m{mM}0\ar[r, "i"']\ar[rrr,"d"',leftrightarrow, bend left]\ar["a"',loop left, leftrightarrow,dotted]& \m{mM}1\ar[d, "i"]\ar[r,"d",leftrightarrow]&\m{AM}2\ar[r,"i"]& \m{AM}3\ar[d,"i"]\ar["a",loop right, dotted] & ~ & ~\\
    ~ & ~ & \m{mM}3\ar[u,"i"]\ar[rrr,"d",leftrightarrow, bend right] & \m{mM}2\ar[l,"i"]\ar[r,"d",leftrightarrow]\ar[r,"a"', shift right=1.5ex, leftrightarrow,dotted]&\m{AM}1\ar[u,"i"] & \m{AM}0\ar[l,"i"] & ~ & ~\\
    \m{MM}1\ar[d,"i"]\ar[rrrruu,"a", leftrightarrow, bend left=40, dotted, crossing over] \ar[d,"d", shift left=1.5ex, leftrightarrow] & \m{MM}0\ar[rrrru,"a",leftrightarrow, bend right, dotted]\ar[l,"i"']\ar[d,"d", shift left=1.5ex, leftrightarrow] & ~ & ~ & ~ & ~ & \m{mm}3\ar[d,"i"]\ar[ldd,"a",dotted, bend right, crossing over] \ar[d,"d",leftrightarrow,shift left=1.5ex] & \m{mm}2\ar[l,"i"]\ar[d,"d",leftrightarrow, shift left=1.5ex]\ar[dddlll,"a",leftrightarrow, dotted, bend left]  \\
    \m{MM}2\ar[r,"i"]\ar[rrruuu,"a",leftrightarrow, dotted, bend left] & \m{MM}3\ar[u,"i"]\ar[ruu,"a",dotted, bend right, crossing over] & ~ & ~ & ~ & ~ & \m{mm}0\ar[r,"i"]\ar[lllld,"a",leftrightarrow,bend right, dotted]& \m{mm}1\ar[u,"i"]\ar[lllldd,"a",leftrightarrow, bend left=40, dotted, crossing over]  \\
    ~ &  ~      &   \m{dm}0\ar[r,"i"]\ar[rrr,"d",leftrightarrow,bend left]  & \m{dm}1\ar[r,"a",leftrightarrow, dotted, shift left=1.5ex]\ar[d,"i"]\ar[r,"d"',leftrightarrow]   & \m{Mm}2\ar[r,"i"] & \m{Mm}3\ar[d,"i"] \\
    ~ & ~       &   \m{dm}3\ar[u,"i"]\ar[rrr,"d",leftrightarrow, bend right]\ar[loop left, "a"',dotted]   & \m{dm}2\ar[l,"i"]\ar[r,"d",leftrightarrow]  & \m{Mm}1\ar[u,"i"] & \m{Mm}0\ar[l,"i"]\ar[loop right, "a"', dotted]
\end{tikzcd}}\end{center}
We see from this diagram that there is a hierarchy of group actions that have different harmonic roles. Not shown in this diagram is the diminished-diminished chord, which is fixed under all these actions. The two connected components of this graph as we have seen correspond to two permutations: $[4,4,3,1]$ for the upper left component and $[4,3,3,2]$ for the lower right component. These two components are isomorphic as graphs, and there is a natural isomorphism which is obtained by the following map:
\begin{align*}
    \m{MM}\longmapsto\m{mm}\\
    \m{mM}\longmapsto\m{Mm}\\
    \m{AM}\longmapsto\m{dm}.
\end{align*}
That is, this is the map that exchanges major with minor and augmented with diminished. If we attempt to extend this map to all chords, we find that the diminished-diminished chord is sent to the augmented-augmented chord, which is just an augmented triad together with the doubling of the root one octave higher. Leaving this aside for the moment, we see that there are four levels of harmonic similarity. At the most basic level, we have inversions. Then major-minor duality maps sets of inversions to sets of inversions. Augmentation is a little different, and highlights the difference between inversion sets. Finally, the isomorphism of the two components corresponding to two permutations provides a mirroring between two very different sets of chords.

Returning to the augmented-diminished operator, we see that it actually highlights the difference between inversions. This is a new facet of four-tone harmonic chords that is not present in three-tone harmonic chords. With three-tone harmonic chords, each inversion of a major or minor triad provides variation for voice leading and some harmonic instability especially for the second inversion. However, with four-tone harmonic chords, the addition of the seventh factor adds truly different harmonic functions to each inversion. For example the root position of a minor-major chord is stable under the augmented-diminished operator, whereas the first inversion of a minor-major chord is sent to the second inversion of a major-major chord. This shows that one inversion can behave much differently than another in music.

\bibliographystyle{plain}
\bibliography{xpapermain.bib}

\addcontentsline{toc}{section}{References}
\end{document}